\documentclass[a4paper, 10pt,oneside]{amsart}
\usepackage{amsmath, amssymb, amsfonts,mathtools}
\usepackage[all]{xy}
\newtheorem{tm}{Theorem}

\newtheorem{kor}[tm]{Corollary}

\newcommand{\card}    {\textrm{card}\,}

\newenvironment{dokaz}
{\noindent\emph{Proof:}\ }
{\hfill $\square$}

\newcommand{\Z}
{{\mathbb Z}}

\newcommand{\C}
{{\mathbb C}}

\newcommand{\g}
{{\mathfrak g}}

\newcommand{\gt}
{\tilde{{\mathfrak g}}}

\newcommand{\hz}
{\hat{\mathfrak h}_\Z}

\newcommand{\h}
{{\mathfrak h}}

\newcommand{\gsl}
{{\mathfrak sl}}

\begin{document}

\author{Goran Trup\v{c}evi\'{c}}
\title[Particle basis of Feigin-Stoyanovsky's type subspaces of level $1$ $\tilde{\gsl}_{\ell+1}(\C)$-modules]{Particle basis of Feigin-Stoyanovsky's type subspaces of level $1$ 
$\tilde{\gsl}_{\ell+1}(\C)$-modules}
\address{Faculty of Teacher Education, University of Zagreb, Trg Matice hrvatske 12, Petrinja, Croatia}
\email{goran.trupcevic@ufzg.hr}

\subjclass[2000]{Primary 17B67; Secondary 05A19.\\ \indent Partially supported by Croatian Science Foundation under the project 2634}
\keywords{affine Lie algebras, principal subspaces, combinatorial basis, character formulas}

\begin{abstract}
We construct particle basis for Feigin-Stoyanovsky's type subspaces of  level $1$ standard $\tilde{\gsl}_{\ell+1}(\C)$-modules. From the description
we obtain character formulas.
\end{abstract}

\maketitle

\section{Introduction}

A problem of finding a monomial basis of a standard module is part
of the Lepowsky-Wilson's program of studying Rogers-Ramanujan type
identities through representation theory of affine Lie algebras
([LW], [LP], [MP]). Description of basis was used to obtain graded dimension
 of these modules, which gave the sum-side
in Rogers-Ramanujan-type identities.

B. Feigin and A. Stoyanovsky initiated another approach to Rogers-Ramanujan type identities 
by considering what they
called a principal subspace of a standard $\tilde{\gsl}_2(\C)$-module
([FS]). These subspaces were further studied  by G. Georgiev ([G]), C. Calinescu, S. Capparelli, 
J. Lepowsky and A. Milas ([CLM1,2], [CalLM1,2,3], [C1,2]), C. Sadowski ([S1,2]), E. Ardonne, 
R. Kedem and M. Stone ([AKS]).

Another type of principal subspace, called Feigin-Stoyanovsky's type subspace, 
was introduced and studied by M. Primc who constructed  a basis of this subspace
and from it he obtained  basis of the whole standard module ([P1,2]). For
$\tilde{\gsl}_{\ell+1}(\C)$, these bases were parameterized
 by $(k,\ell+1)$-admissible configurations ([FJLMM]), combinatorial objects 
introduced
and further studied by Feigin et al. in [FJLMM] and [FJMMT], where bosonic and fermionic formulas
for characters were obtained. Primc and M. Jerkovic obtained fermionic formulas for characters of standard 
$\tilde{\gsl}_3(\C)$-modules by using intertwining operators and admissible configurations ([J2]), or by
using quasi-particle bases [JP]). In our previous work ([T]), we have used $(1,\ell+1)$-admissible configurations
to combinatorially obtain character formulas for Feigin-Stoyanovsky's type subspaces of level $1$ standard 
$\tilde{\gsl}_{\ell+1}(\C)$-modules. In this note we use an approach similar to Georgiev and Jerkovi\' c and Primc to 
construct particle basis for Feigin-Stoyanovsky's type subspaces of level $1$ standard $\tilde{\gsl}_{\ell+1}(\C)$-modules. From this description
we immediately obtain character formulas.

\section{Affine Lie algebra $\hat{\mathfrak sl}_{\ell+1} (\C)$}

\label{FS_sect}

Let ${\mathfrak g}={\mathfrak sl}_{\ell+1} (\C)$ be a simple finite-dimensional Lie algebra of  type $A_\ell$. 
Fix a Cartan subalgebra $\h\subset\g$ and denote by $R$ the corresponding root system; $\g$ has a
root decomposition $\g=\h \oplus \coprod_{\alpha\in R} \g_\alpha$. 
Let $\Pi=\{\alpha_1,\dots,\alpha_\ell\}$ 
be a basis of the root system $R$, and let $\{\omega_1,\dots,\omega_\ell\}$ be the
corresponding set of fundamental weights, $\langle\omega_i,\alpha_j\rangle=\delta_{ij}$.
Set $\omega_0=0$ for convenience. Let  $\langle
\cdot,\cdot\rangle$ be a a normalized invariant bilinear form on
$\g$; we identify $\h$ with $\h^*$ via $\langle
\cdot,\cdot\rangle$.
Denote by $Q$
the root lattice, and by $P$ the weight lattice
 of $\g$.  Also for each root $\alpha\in R$ fix a root vector $x_\alpha\in\g_\alpha$.

Let $\gt$ be the associated untwisted affine Lie algebra ([K]),
$$\gt=\g\otimes \C[t,t^{-1}]\oplus \C c \oplus \C d.$$
Denote by $x(m)=x\otimes t^m$ for $x\in\g$, $m\in\Z$, and define formal Laurent series 
$x(z)=\sum_{m\in \Z} x(m)z^{-m-1}$.
Denote by  $\Lambda_0$, \dots,
$\Lambda_\ell$ fundamental weights for $\gt$.

Fix a minuscule weight $\omega = \omega_\ell$ and set 
$$\Gamma = \{\,\alpha \in R \mid \langle\alpha,\omega\rangle = 1\}=\{\alpha_i+\dots+\alpha_\ell \mid i=1,\dots,\ell\}.$$ 
Denote by  $\gamma_i=\alpha_i+\dots+\alpha_\ell$. 
Then
$$
\mathfrak g  =
\mathfrak g_{-1} \oplus \mathfrak g_0 \oplus \mathfrak g_1, \qquad {\mathfrak g}_0  =   {\mathfrak h} \oplus
\sum_{\langle \alpha, \omega \rangle =0}\, {\mathfrak g}_\alpha,\quad
 {\mathfrak g}_{\pm1} =
\sum_{\alpha \in \pm \Gamma}\, {\mathfrak g}_\alpha,$$
 is a $\mathbb Z$-gradation of ${\mathfrak g}$.
Subalgebras ${\mathfrak g}_1$ and ${\mathfrak g}_{-1}$ are
commutative. 

 The $\mathbb Z$-gradation of ${\mathfrak g}$ gives the $\Z$-gradation of the affine Lie
algebra $\gt$:
$$\nonumber  \gt = \gt_{-1} + \gt_0 + \gt_1,\qquad  \gt_0 = {\mathfrak g}_0\otimes\C [t,t^{-1}]\oplus \C c \oplus \C d,\quad
\gt_{\pm 1} = {\mathfrak g}_{\pm 1}\otimes\C [t,t^{-1}]. $$
Again, $\gt_{-1}$ and $\gt_1$ are commutative subalgebras. 
We will call elements  $\gamma\in\Gamma$ {\em colors} and we will say that $x_\gamma(-m)$ is an element of {\em color} $\gamma$ and {\em degree} $m$.

Let $L(\Lambda_r)$ be a standard (i.e. integrable highest weight) $\gt$-module
of level $1$. Denote by $v_r$ the highest weight vector of $L(\Lambda_r)$.
Define a
\emph{Feigin-Stoyanovsky's type subspace}
$$W(\Lambda_r)=U(\gt_1)\cdot v_r\subset L(\Lambda_r).$$
By Poincar\'e-Birkhoff-Witt theorem, we have a spanning set of
$W(\Lambda_r)$ consisting of monomial vectors
\begin{equation} \label{MVect_def}
\{b v_r \mid b=x_\ell(-m_{\ell,n_\ell})\cdots x_1(-m_{1,n_1})\cdots x_1(-m_{1,1}), m_{i,j+1}\geq m_{i,j}>0, n_i\geq 0 \},
\end{equation}
where we write $x_i$ instead of $x_{\gamma_i}$, for short.

\section{VOA construction}

We briefly recall the vertex operator algebra construction
of standard $\gt$-modules $L(\Lambda_r)$ from [FK], [S]. For details and notation we turn the reader to [FLM], [DL] and [LL].

Consider tensor products
$V_P  =  M(1)\otimes \C [P]$ and $V_Q  =  M(1)\otimes \C [Q]$,
where $M(1)$ is the Fock space for the Heisenberg subalgebra
$\hz=\sum_{n\in\Z \setminus\{0\}}\h\otimes t^n \oplus \C c$, and $\C
[P]$ and $\C [Q]$ are group algebras of the weight and root lattice
with bases consisting of $\{e^\lambda\,|\,\lambda\in P\}$, and
$\{e^\alpha\,|\,\alpha\in Q\}$, respectively. We identify $\C [P]$ with $1\otimes \C [P]\subset V_P$.

Space $V_Q$ has a structure of vertex operator algebra and
$V_P$ is a module for this algebra:
\begin{equation}\label{VOA_rel}
Y(e^\lambda,z)=E^-(-\lambda,z)E^+(-\lambda,z)\otimes
e^\lambda z^\lambda \epsilon(\lambda,\cdot),
\end{equation} 
where $ E^{\pm}(\lambda,z)=\exp \left(\sum_{m\geq 1}\lambda(\pm m) {z^{\mp m}}/{\pm m}\right)$,
$e^\lambda$
is a multiplication operator, $z^\lambda \cdot e^\mu=e^\mu z^{\langle \lambda,\mu \rangle}$ and $\epsilon(\cdot,\cdot)$ is a
$2$-cocycle (cf. [FLM]).

By using vertex operators, one can define the structure of
$\gt$-module on $V_P$ by setting $x_\alpha(z)=Y(e^\alpha,z)$ for $\alpha\in R$.
This gives $V_Q\simeq L(\Lambda_0)$ and $V_Q e^{\omega_r}\simeq L(\Lambda_r)$, with highest weight vectors
$v_0=1$ and $v_r=e^{\omega_r}$, and $V_P\cong L(\Lambda_0)\oplus\dots
\oplus L(\Lambda_\ell)$.

From vertex operator formula \eqref{VOA_rel} one  easily obtains  the following
relations on $L(\Lambda_r)$
\begin{eqnarray}
\label{DC_rel} x_i^2(z)& = & 0, \quad 1\leq i \leq \ell, \\
\label{Interaction_rel} x_i(z) x_j(z) & = & 0, \quad 1\leq i < j \leq \ell,\\
\label{Init_rel} x_i(m)v_r & = & 0, \quad m  \geq -\delta_{i\leq r},\\
\label{Init2_rel} x_r(-1) v_{r-1}& = & C e^{\omega_\ell+\omega_r}=C e^{\omega_\ell} v_r,
\end{eqnarray}
for some $C\in\C^\times$. Here, $\delta_{i\leq j}$ is $1$ if $i\leq j$, $0$ otherwise.

For the proof of linear independence we will be using certain coefficients of intertwining operators
$${\mathcal Y} (e^\lambda,z)=Y(e^\lambda,z)e^{i\pi\lambda}c(\cdot,\lambda),$$
for $\lambda\in P$, where $c(\cdot,\lambda)$ is a commutator map
(cf. [DL]). Let $\lambda_i=\omega_i-\omega_{i-1}$ for $i=1,\dots,\ell$.
From Jacobi identity ([DL]) we see that operators ${\mathcal Y} (e^{\lambda_i},z)$ commute with the action of
$\gt_1$. Define the following coefficients of intertwining operators (cf. [P3])
$$[i]=\textrm{Res}\, z^{-1-\langle \lambda _i, \omega_{i-1} \rangle} {\mathcal Y} (e^{\lambda_i},z),
$$ for $i=1,\dots,\ell$.
From \eqref{VOA_rel}, it follows (cf. [P3])
\begin{equation}\label{Inttw_rel}
[i] v_{i-1}=Cv_i,
\end{equation}
for some $C\in\C^\times$.

We will also be using simple current operators $e^{\omega_i}$, $i=1,\dots,\ell$.
For $\alpha\in R$ and $\lambda \in P$ from \eqref{VOA_rel} we get the following commutation relation
$$x_\alpha(z) e^\lambda=\epsilon(\alpha, \lambda)z^{\langle \alpha,\lambda \rangle} e^\lambda x_\alpha(z).$$
By comparing coefficients, we get 
$$x_\alpha(m) e^\lambda=\epsilon(\alpha, \lambda) e^\lambda x_\alpha(m+\langle \alpha,\lambda \rangle).$$
In particular, for $\alpha=\gamma_i$ and $\lambda=\omega_j$, we get
\begin{equation} \label{komut_ea_xn_jed}
x_i(m) e^{\omega_j}=\epsilon(\gamma_i, \omega_j) e^{\omega_j} x_i(m+ \delta_{i\leq j}).
\end{equation}
%
%

\section{Basis of $W(\Lambda_r)$}

To reduce the spanning set \eqref{MVect_def} and to prove linear independence, we need a linear order on monomials.
Define a linear order $x_i(n) < x_j(m)$ if either $i>j$ or $i=j$ and $n<m$. We assume that in all monomials factors are 
sorted descendingly from right to left, like in \eqref{MVect_def}. We compare two monomials $b_1$ and $b_2$ 
by comparing their factors from right to left (reverse lexicographic order): $b_1<b_2$ if either $b_2=bb_1$ or 
$b_1=b_1'x_i (n) b$, $b_2=b_2'x_j (m) b$ and $x_i(n) < x_j(m)$, for some monomials $b,b_1,b_2$. This linear order is compatible 
with multiplication: if $b>c$, then $ab>ac$.

For a monomial $b=x_\ell(-m_{\ell,n_\ell})\cdots x_\ell(-m_{\ell,1})\cdots x_1(-m_{1,n_1})\cdots x_1(-m_{1,1})$ define its {\em degree}, {\em weight}
and {\em length} by
$d(b)=m_{\ell,n_\ell}+\dots+m_{\ell,1}+\dots+m_{1,n_1} +\dots + m_{1,1}$,  $w(b)=n_1 \gamma_1+\dots +n_\ell \gamma_\ell$ and
$l(b)=n_1 +\dots +n_\ell$.

\begin{tm} \label{Span1_tm}
A spanning set of $W(\Lambda_r)$ is given by the set of monomial vectors \eqref{MVect_def}
satisfying initial conditions
\begin{equation} \label{IC_ineq}
 m_{i,n}  \geq  1 +\sum_{i<j} n_i + \delta_{j\leq r}
\end{equation}
and difference conditions
\begin{equation} \label{DC_ineq}
 m_{i,n+1}  \geq  m_{i,n}+2,\quad 1\leq n \leq n_i -1.
\end{equation}
\end{tm}
\begin{dokaz}
Difference conditions follow from \eqref{DC_rel}: Assume that $b$ doesn't satisfy \eqref{DC_ineq}. Then 
$b=b' x_j(-m)x_j(-m')$, for some monomial $b'$ and $m' \leq m \leq m' + 1$. By \eqref{DC_rel} and \eqref{Init_rel}, on $W(\Lambda_r)$ we have
$$x_j(-m)x_j(-m') = C_1 x_j(-m-1)x_j(-m'+1)+ \dots + C_{m'-1}  x_j(-m-m'+1) x_j(-1), $$
for some $C_i\in \C^\times $. 
Multiply this by $b'$ to obtain $b$ expressed as a linear combination of greater monomials of the same degree and weight.

Now assume that $b$ doesn't satisfy \eqref{IC_ineq}; let $b = b_2 x_j(-m) b_1$ where 
$b_1$ contains all factors of colors $\gamma_1,\dots,\gamma_{j-1}$ and
$$ m <  1 +\sum_{i<j} n_i + \delta_{j\leq r}. $$
We will prove that $b$ can be expressed in terms of greater monomials of the same degree and weight. 
The  proof is done by induction on the length $l(b_1)=\sum_{i<j} n_i$.
If $l(b_1) =0$, then \eqref{Init_rel} gives $x_j(-m)v_r=0$.
Now, assume that all monomials $a$ with $l(a_1)<l(b_1)$ can be 
expressed with greater monomials of the same degree and weight.   
We can also assume that $ m = \sum_{i<j} n_i + \delta_{j\leq r}$. Let $x_k(-n)$ be the smallest factor in $b_1$;
$b_1= x_k(-n) b_1'$. 
By \eqref{Interaction_rel} and \eqref{Init_rel} we   have
\begin{equation}
\begin{multlined}
x_j(-m) x_k(-n)=  C_{1+ \delta_{j\leq r}} x_j(-1- \delta_{j\leq r}) x_k(-n-m +1+\delta_{j\leq r} )+ \dots \\ 
 + C_{m-1} x_j(-m+1)   x_k(-n-1 )+ C_{m+1} x_j(-m-1)x_k(-n +1 )+\dots
\end{multlined}
\end{equation}
for some $C_i\in \C $. 
Multiply this with $b_2 b_1'$ and obtain
\begin{equation}
\begin{multlined}
b = b_2 x_j(-1- \delta_{j\leq r}) x_k(-n-m +1+\delta_{j\leq r} ) b_1'+ \dots \\ 
 + b_2 x_j(-m+1)   x_k(-n-1 ) b_1' + b_2 x_j(-m-1)x_k(-n +1 ) b_1'+\dots
\end{multlined}
\end{equation}
On the right-hand side we have monomials
$$b_2 x_j(-m-1)x_k(-n +1 ) b_1', b_2 x_j(-m-2)x_k(-n +2 ) b_1',\dots$$ 
which are greater 
than $b$. But we also have the first few monomials 
\begin{equation}\label{Prvi_mon_list_2}
b_2 x_j(-1- \delta_{j\leq r}) x_k(-n-m +1+\delta_{j\leq r} ) b_1', \dots,  
 b_2 x_j(-m+1) x_k(-n-1 ) b_1'.
\end{equation}
Consider their factors $x_j(-1- \delta_{j\leq r}) b_1', \dots, x_j(-m+1) b_1'$.
By the inductive assumption, they can be expressed as linear 
combinations of greater monomials of the same degree and weight. 
Then it is obvious that by multiplying these linear expressions by 
$b_2  x_k(-n-m +1+\delta_{j\leq r} )$, \dots, $b_2 x_k(-n-1)$ we obtain linear 
expressions for monomials in \eqref{Prvi_mon_list_2} in terms of greater monomials. 
Moreover, these monomials will also be greater than $b$.
Hence, we have expressed $b$ in terms of greater monomial.
\end{dokaz}

\begin{tm} 
A spanning set
\begin{equation} \label{basis_def}
\mathcal B=\{ b v_r | b\ \textrm{satisfies \eqref{DC_ineq} and \eqref{IC_ineq}}\}\end{equation}
is a basis of $W(\Lambda_r)$.
\end{tm}

\begin{dokaz}
Let $b\in \mathcal B$. We first prove a particular case: if
$$Cbv_r=0$$
then $C=0$. We prove this by induction on degree of $b$.
Let $x_i(-n)$ be the greatest factor in $b$; $b=b' x_i(-n)$.

If $i\leq r$ then, since $v_r=e^{\omega_r}=e^{\omega_r}v_0$, we have
$$Cbv_r=Cbe^{\omega_r}v_0=e^{\omega_r}C'Cb''v_0=0,$$
where $C'\in\C^\times$ and $b''$ is obtained from $b$ by decreasing degrees of factors of color $\gamma_1,\dots,\gamma_r$
by $1$ (see \eqref{komut_ea_xn_jed}). Since $e^{\omega_r}$ is injective, $Cb''v_0=0$. Monomial $b''$ satisfies difference and initial conditions
for $W(\Lambda_0)$ and it is of a smaller degree  then $b$. By the inductive assumption, we conclude that $C=0$.

If $i>r$ and $n>1$ then use operators $[i][i-1]\cdots [r+1]$ to obtain
$$Cbv_{i}=0.$$ Then, by \eqref{komut_ea_xn_jed},
$$Cbv_{i}= Cbe^{\omega_i}v_0=e^{\omega_i}C'Cb''v_{0}=0,$$
where $C'\in\C^\times$ and  $b''$ is obtained from $b$ by decreasing degrees of factors of color $\gamma_i$
by $1$ (see \eqref{Inttw_rel} and \eqref{komut_ea_xn_jed}). Again, $b''$ satisfies difference and initial conditions
for $W(\Lambda_0)$ and it is of a smaller degree   then $b$. By induction, we conclude that $C=0$.

If $i>r$ and $n=1$ then use operators $[i-1]\cdots [r+1]$ to obtain
$$Cbv_{i-1}=0$$ 
(see \eqref{Inttw_rel}). By \eqref{Init2_rel} and \eqref{komut_ea_xn_jed}, we have
$$C b' x_i(-1)v_i=C' C b'e^{\omega_\ell}v_i=e^{\omega_\ell} C'' Cb''v_{i}=0,$$
where  $C',C''\in\C^\times$ and  $b''$ is obtained from $b'$ by decreasing degrees of all factors
by $1$. Monomial $b''$ satisfies difference and initial conditions
for $W(\Lambda_i)$ and and it is of a smaller degree   then $b$. By induction, we conclude that $C=0$.

Now turn to the general relation of linear dependence:
\begin{equation} \label{LD_rel}
\sum_b C_b bv_r=0.
\end{equation}
Let $b_{min}$ be the smallest monimial in \eqref{LD_rel}. We use the same operators as in the previous case
to peel down $C_{b_{min}} b_{min} v_r$ to $C_{b_{min}} v_j$, for some $j$. Note that operators that we used above 
at some point annihilate other monomial vectors in \eqref{LD_rel}. We will get 
$$C_{b_{min}} v_j=0$$
and we conclude $C_{b_{min}}=0$. We proceed inductively to conclude that all coefficients $C_b$ in \eqref{LD_rel}
are $0$.   
\end{dokaz}

From this combinatorial description of basis of $W(\Lambda_r)$ we immediately obtain character formulas (cf. [J1],
[T]):  for $n_1, \dots , n_\ell\geq 0$, set $\alpha=n_ 1 \gamma_1+\dots+n_\ell \gamma_\ell$ and 
$\chi_{W(\Lambda_r)}^\alpha(q)=\sum_i q^i \card \{b \mid w(b)=\alpha, d(b)=i \}$
\begin{kor} 
$$
\chi_{W(\Lambda_r)}^\alpha(q)=\frac{ q^{\sum_{i=1}^\ell n_i^2+\sum_{1\leq i < j \leq \ell} n_i n_j+\sum_{i=1}^r n_i}} {(q)_{n_1}(q)_{n_2}\cdots(q)_{n_\ell}},
$$
where $(q)_n=(1-q)\cdots(1-q^n)$.
 \end{kor}

\end{document}